\overfullrule=0pt
\centerline {\bf Minimax theorems in a fully non-convex setting}\par
\bigskip
\bigskip
\centerline {\it Dedicated to Professor Wataru Takahashi, with esteem and friendship, on his 75th birthday}
\bigskip
\bigskip
\centerline {BIAGIO RICCERI}\par
\bigskip
\bigskip
{\bf Abstract.} In this paper, we establish two minimax theorems for functions $f:X\times I\to {\bf R}$, where $I$ is a real interval, without
assuming that $f(x,\cdot)$ is quasi-concave. Also, some related applications are presented.\par
\bigskip
\bigskip
The most known minimax theorem ([7]) ensures the occurrence of the equality
$$\sup_Y\inf_Xf=\inf_X\sup_Yf$$
for a function $f:X\times Y\to {\bf R}$ under the following assumptions: $X$, $Y$
are convex sets in Hausdorff topological vector spaces, one of them is compact,
$f$ is lower semicontinuous and quasi-convex in $X$, and upper semicontinuous and
quasi-concave in $Y$.\par
\smallskip
In the past years, we provided some contributions to the subject where, keeping the assumption of quasi-concavity on
$f(x,\cdot)$, we proposed alternative hypotheses on $f(\cdot,y)$. Precisely, in [2],
we assumed the inf-connectedness of $f(\cdot,y)$ and, the same time, that $Y$ is a real interval, while,
in [5], we assumed the inf-compactness and uniqueness of the global minimum of $f(\cdot,y)$.\par
\smallskip
In the present paper, we offer a new contribution where the hypothesis that $f(x,\cdot)$ is
quasi-concave is no longer assumed.\par
\smallskip
Let $T$ be a topological space. A function $g:T\to [-\infty,+\infty[$ is said to be relatively inf-compact if, for each $r\in {\bf R}$,
there exists a compact set $K\subseteq T$ such that $g^{-1}(]-\infty,r[)\subseteq K$. Moreover, $g$ is said to be
inf-connected if, for each $r\in {\bf R}$, the set $g^{-1}(]-\infty,r[)$ is connected. For the basic notions on multifunctions,
we refer to [1].\par
\smallskip
Our main results are as follows:\par
\medskip
THEOREM 1. - {\it Let $X$ be a topological space, let $I$ be a real interval and let $f:X\times  I\to {\bf R}$ be a continuous function such
that, for each $\lambda\in I$, the set of all global minima of the function $f(\cdot,\lambda)$ is connected. Moreover, assume that
there exists a non-decreasing sequence of compact intervals, $\{I_n\}$, with $I=\cup_{n\in {\bf N}}I_n$, such that, for each $n\in {\bf N}$,
the following conditions are satisfied:\par
\noindent
$(a_1)$\hskip 5pt the function $\inf_{\lambda\in I_n}f(\cdot,\lambda)$ is relatively inf-compact\ ;\par
\noindent
$(b_1)$\hskip 5pt for each $x\in X$, the set of all global maxima of the restriction of the function $f(x,\cdot)$ to $I_n$ is connected.\par
Then, one has
$$\sup_Y\inf_X=\inf_X\sup_Yf\ .$$}\par
\medskip
THEOREM 2. - {\it Let $X$ be a topological space, let $I$ be a compact real interval and let $f:X\times I\to {\bf R}$ be an upper semicontinuous
 function such that $f(\cdot, \lambda)$ is continuous for all $\lambda\in I$.
Assume that:\par
\noindent
$(a_2)$\hskip 5pt there exists a set $D\subseteq I$, dense in $I$,  such that the function $f(\cdot,\lambda)$ is inf-connected for
all $\lambda\in D$\ ;\par
\noindent
$(b_2)$\hskip 5pt for each $x\in X$, the set of all global maxima of the function $f(x,\cdot)$ is connected.\par
Then, one has
$$\sup_Y\inf_X=\inf_X\sup_Yf\ .$$}\par
\medskip
REMARK 1. -
We want to remark that, in both Theorems 1 and 2, it is essential that $I$ be a real interval. To see this, consider the following
example. Take
$$X=I=\{(t,s)\in {\bf R}^2 : t^2+s^2=1\}$$
and define $f:X\times I\to {\bf R}$ by
$$f(t,s,u,v)=tu+sv$$
for all $(t,s), (u,v)\in X$. Clearly, $f$ is continuous, $f(\cdot,\cdot,u,v)$ is inf-connected and has a unique global minimum, and
$f(t,s,\cdot,\cdot)$ has a unique global maximum. However, we have
$$\sup_X\inf_If=-1<1=\inf_X\sup_If\ .$$
\medskip
The common key tool in our proofs of Theorems 1 and 2 is provided by the following general principle:\par
\medskip 
THEOREM A ([2], Theorem 2.2). - {\it Let $X$ be a topological space, let $I$ be a compact real interval and let
$S\subseteq X\times I$ be a connected set whose projection on $I$ is the whole of $I$.\par
Then, for every upper semicontinuous multifunction $\Phi:X\to 2^I$, with non-empty, closed and connected values, the graph of $\Phi$ intersects $S$.}\par
\medskip
Another known proposition which is used in the proof of Theorem 1 is as follows:\par
\medskip
PROPOSITION A ([5], Proposition 2.1). - {\it Let $X$ be a topological space, $Y$ a non-empty set, $y_0\in Y$ and $f:X\times Y\to {\bf R}$
a function such that $f(\cdot,y)$ is lower semicontinuous  for
all $y\in Y$ and relatively inf-compact  for $y=y_0$. Assume also that there is a non-decreasing sequence of sets $\{Y_n\}$,
with $Y=\cup_{n\in {\bf N}}Y_n$,
such that
$$\sup_{Y_n}\inf_Xf=\inf_X\sup_{Y_n}f$$
for all $n\in {\bf N}$.\par
Then, one has
$$\sup_Y\inf_Xf=\inf_X\sup_Yf\ .$$}\par
\medskip
A further result which is used in the proofs of Theorems 1 and 2 is provided by the following proposition which, in
the given generality, is new:\par
\medskip
PROPOSITION 1. - {\it Let $X, Y$ be two topological spaces and let $f:X\times Y\to {\bf R}$ be
a lower semicontinuous function such that $f(x,\cdot)$ is continuous for all $x\in X$.
Moreover, assume that, for each $y\in Y$, there exists
a neighbourhood $V$ of $y$ such that the function $\inf_{v\in V}f(\cdot,v)$ is
relatively inf-compact. For each $y\in Y$, set
$$F(y)=\left \{u\in X : f(u,y)=\inf_{x\in X}f(x,y)\right \}\ .$$
Then, the multifunction $F$ is upper semicontinuous.}\par
\smallskip
PROOF. Let $C\subseteq X$ be a closed set. We have to prove that $F^-(C)$ is closed. So, let $\{y_{\alpha}\}_{\alpha\in D}$ be a net in $F^-(C)$ converging
to some $\tilde y\in Y$. For each $\alpha\in D$, pick $u_{\alpha}\in F(y_{\alpha})\cap C$.  By assumption,
there is a neighbourhood $V$ of $\tilde y$ such that the function $\inf_{v\in V}f(\cdot,v)$ is relatively
inf-compact. Since the function $\inf_{x\in X}f(x,\cdot)$ is upper semicontinuous, we can assume that it is bounded above on $V$. Fix $\rho>\sup_V\inf_Xf$.
 Then, there is a compact set $K\subseteq X$ such that
$$\left \{x\in X : \inf_{v\in V}f(x,v)<\rho\right \}\subseteq K\ .$$
But 
$$\left \{x\in X : \inf_{v\in V}f(x,v)<\rho\right \}=\bigcup_{v\in V}\left \{x\in X : f(x,v)<\rho\right\}$$
and so
$$\bigcup_{v\in V}\left \{x\in X : f(x,v)<\rho\right\}\subseteq K\ .\eqno{(1)}$$
Let $\alpha_1\in D$ be such that $y_{\alpha}\in V$ for all $\alpha\geq \alpha_1$. Consequently, by $(1)$,
$u_{\alpha}\in K$ for all $\alpha\geq\alpha_1$. By compactness, the net $\{u_{\alpha}\}_{\alpha\in D}$ has a cluster point $\tilde u\in K$.
Clearly, $(\tilde u,\tilde y)$ is a cluster point in $X\times Y$ of the net $\{(u_{\alpha},y_{\alpha})\}_{\alpha\in D}$. We claim that
$$f(\tilde u,\tilde y)\leq\limsup_{\alpha}f(u_{\alpha},y_{\alpha})\ .$$
Arguing by contradiction, assume the contrary and fix $r$ so that
$$\limsup_{\alpha}f(u_{\alpha},y_{\alpha})<r<f(\tilde u,\tilde y)\ .$$
Then,  there would be $\alpha_2\in D$ such that 
$$f(u_{\alpha},y_{\alpha})<r$$
for all $\alpha\geq\alpha_2$. On the other hand, since, by assumption, the set $f^{-1}(]r,+\infty[)$ is open, there would be
$\alpha_3\geq\alpha_2$ such that
$$r<f(u_{\alpha_3},y_{\alpha_3})$$
which gives a contradiction. Now, fix $x\in X$. Then, since $u_{\alpha}\in F(y_{\alpha})$, we have
$$f(\tilde u,\tilde y)\leq \limsup_{\alpha}f(u_{\alpha},y_{\alpha})\leq \lim_{\alpha}f(x,y_{\alpha})=f(x,\tilde y)\ .$$
That is, $\tilde u\in F(\tilde y)$. Since $C$ is closed, $\tilde u\in C$. Hence, $\tilde y\in F^-(C)$ and this ends the proof.\hfill $\bigtriangleup$\par
\medskip
We now can prove Theorems 1 and 2.\par
\medskip
{\it Proof of Theorem 1.} Fix $n\in {\bf N}$. Let us prove that
$$\sup_{I_n}\inf_Xf=\inf_X\sup_{I_n}f\ .\eqno{(2)}$$
Consider the multifunction $F:I_n\to 2^X$ defined by
$$F(\lambda)=\left \{u\in X : f(u,\lambda)=\inf_{x\in X}f(x,\lambda)\right\}$$
for all $\lambda\in I_n$. Thanks to Proposition 1, $F$ is upper semicontinuous and, by assumption, its values are non-empty, compact and
connected. As a consequence, by Theorem 7.4.4 of [1], the graph of $F$ is connected. Let $S$ denote the graph of the inverse of $F$. So,
$S$ is connected as it is homeomorphic to the graph of $F$. Now, consider the multifunction $\Phi:X\to 2^{I_n}$ defined by
$$\Phi(x)=\left\{\mu\in I_n : f(x,\mu)=\sup_{\lambda\in I_n}f(x,\lambda)\right\}$$
for all $x\in X$. By Proposition 1 again, the multifunction $\Phi$ is upper semicontinuous and, by assumption, its values are non-empty, closed and
connected. After noticing that the projection of $S$ on $I_n$ is the whole of $I_n$, we can apply Theorem A. Therefore, there exists $(\tilde x,
\tilde\lambda)\in S$ such that $\tilde\lambda\in \Phi(\tilde x)$. That is
$$f(\tilde x,\tilde\lambda)=\inf_{x\in X}f(x,\tilde\lambda)=\sup_{\lambda\in I_n}f(\tilde x,\lambda)\ .\eqno{(3)}$$
Clearly, $(2)$ follows from $(3)$. Now, the conclusion is a direct consequence of Proposition A.\hfill $\bigtriangleup$
\medskip
{\it Proof of Theorem 2}. Arguing by contradiction, assume the contrary and fix a constant $r$ so that
$$\sup_I\inf_Xf<r<\inf_X\sup_If\ .$$
Let $G:I\to 2^X$ be the multifunction defined by
$$G(\lambda)=\{x\in X : f(x,\lambda)<r\}$$
for all $\lambda\in I$. Notice that $G(\lambda)$ is non-empty for all $\lambda\in I$ and connected for all $\lambda\in D$. Moreover, the
graph of $G$ is open in $X\times I$ and so $G$ is lower semicontinuous. Then, by Proposition 5.7 of [3], the graph of $G$ is connected and so
the graph of the inverse of $G$, say $S$, is connected too. Consider the multifunction $\Phi:X\to 2^I$ defined by
$$\Phi(x)=\left\{\mu\in I : f(x,\mu)=\sup_{\lambda\in I}f(x,\lambda)\right\}$$
for all $x\in X$. Notice that $\Phi(x)$ is non-empty, closed and connected, in view of $(b_2)$. By Proposition 1, the multifunction $\Phi$ is upper semicontinuous.
Now, we can apply Theorem A. So, there exists $(\hat x,\hat\lambda)\in S$ such that $\hat\lambda\in \Phi(\hat x)$. This implies that
$$f(\hat x,\hat\lambda)<r<\inf_X\sup_If\leq \sup_{\lambda\in I}f(\hat x,\lambda)=f(\hat x,\hat\lambda)$$
which is absurd.\hfill $\bigtriangleup$
\medskip
Here is an application of Theorem 1.\par
\medskip
THEOREM 3. - {\it Let $(H,\langle\cdot,\cdot\rangle)$ be a real inner product space, let $K\subset H$ be a compact and convex set, 
with $0\notin K$, and
let $f:X\to K$ be a continuous function, where
$$X=\bigcup_{\lambda\in {\bf R}}\lambda K\ .$$
Assume that there are two numbers $\alpha, c$, with 
$$\inf_{x\in X}\|f(x)\|<c<\|f(0)\|\ ,$$
such that:\par
\noindent
$(a)$\hskip 5pt $\{x\in X : \langle x,f(x)\rangle=\alpha\}\subset \{x\in X : \|f(x)\|<c\}$\ ;\par
\noindent
$(b)$\hskip 5pt $\{x\in X :  c^2\langle x,f(x)\rangle=\alpha\|f(x)\|^2\}\subset \{x\in X : \|f(x)\|\geq c\}\ .$\par
Then, there exists $\tilde\lambda\in {\bf R}$ such that the set
$$\{x\in X : x=\tilde\lambda f(x)\}$$
is disconnected.}\par
\smallskip
PROOF. Consider the function $\varphi:X\times {\bf R}\to {\bf R}$ defined by
$$\varphi(x,\lambda)=\|x-\lambda f(x)\|^2-c^2\lambda^2+2\alpha\lambda$$
for all $(x,\lambda)\in X\times {\bf R}$. 
Notice that
$$\varphi(x,\lambda)=\|x\|^2+(\|f(x)\|^2-c^2)\lambda^2-2(\langle x,f(x)\rangle-\alpha)\lambda\ .\eqno{(4)}$$
Further, observe that, when $\|f(x)\|\geq c$, in view of $(a)$, we have
$$\sup_{\lambda\in {\bf R}}\varphi(x,\lambda)=+\infty\eqno{(5)}$$
as well as
$$\varphi(x,-\lambda)\neq \varphi(x,\lambda)\eqno{(6)}$$
for all $\lambda>0$. When $\|f(x)\|\geq c$ again, the function $\varphi(x,\cdot)$ is convex and so, by $(6)$, for each $\lambda>0$,
its restriction to $[-\lambda,\lambda]$ it has a unique global maximum. Clearly, $\varphi(x,\cdot)$ has
the same uniqueness property also when $\|f(x)\|<c$.
Now, observe that, for each $\lambda\in {\bf R}$, the function $\lambda f$ has a fixed point in $X$, in view of the Schauder theorem. Hence, we have
$$\sup_{\lambda\in {\bf R}}\inf_{x\in X}\varphi(x,\lambda)=\sup_{\lambda\in {\bf R}}(-c^2\lambda^2+2\alpha\lambda)={{\alpha^2}\over {c^2}}\ .\eqno{(7)}$$
We claim that
$${{\alpha^2}\over {c^2}}<\inf_{x\in X}\sup_{\lambda\in {\bf R}}\varphi(x,\lambda)\ .\eqno{(8)}$$
First, observe that, since $0\notin K$, every closed and bounded subset of $X$ is compact. This easily implies that, for each $\mu>0$,
the function $x\to \inf_{|\lambda|\leq \mu}\varphi(x,\lambda)$ is relatively inf-compact.
Consequently,
the sublevel sets of the function $x\to\sup_{\lambda\in {\bf R}}\varphi(x,\lambda)$ (which is finite
if $\|f(x)\|<c$) are compact. Therefore,
there exists $\tilde x\in X$ such that
$$\sup_{\lambda\in {\bf R}}\varphi(\tilde x,\lambda)=\inf_{x\in X}\sup_{\lambda\in {\bf R}}\varphi(x,\lambda)\ .\eqno{(9)}$$
So, by $(5)$, one has $\|f(\tilde x)\|<c$.
Clearly, we also have
$$\sup_{\lambda\in {\bf R}}\varphi(\tilde x,\lambda)=\|\tilde x\|^2+{{|\langle \tilde x,f(\tilde x)\rangle-\alpha|^2}\over {c^2-\|f(\tilde x)\|^2}}\ .\eqno{(10)}$$
Let us prove that
$$\|\tilde x\|^2+{{|\langle \tilde x,f(\tilde x)\rangle-\alpha|^2}\over {c^2-\|f(\tilde x)\|^2}}>{{\alpha^2}\over {c^2}}\ .\eqno{(11)}$$
After some manipulations, one realizes that $(11)$ is equivalent to
$${{1}\over {c^2-\|f(\tilde x)\|^2}}\left (2\alpha\langle \tilde x,f(\tilde x)\rangle-|\langle \tilde x,f(\tilde x)\rangle|^2-{{\alpha^2}\over {c^2}}\|f(\tilde x)\|^2\right )<\|\tilde x\|^2\ .\eqno{(12)}$$
Now, for each $y\in X\setminus \{0\}$, $t\in {\bf R}$, set
$$I(y,t)=\{x\in H : \langle x,y\rangle=t\}\ .$$
Consider the inequality
$${{1}\over {c^2-\|y\|^2}}\left ( 2\alpha t-t^2-{{\alpha^2}\over {c^2}}\|y\|^2\right )<{{t^2}\over {\|y\|^2}}\ .\eqno{(13)}$$
After some manipulations, one realizes that $(13)$ is equivalent to
$$(\alpha\|y\|^2-tc^2)^2>0\ .$$
So, $(13)$ is satisfied if and only if
$$\alpha\|y\|^2\neq tc^2\ .\eqno{(14)}$$
Observe that
$${{|t|}\over {\|y\|}}=\hbox {\rm dist}(0,I(y,t))\leq \hbox {\rm dist}(0,I(y,t)\cap X)\ .\eqno{(15)}$$
Therefore, if $(14)$ is satisfied, for each $x\in I(y,t)\cap X$, in view of $(13)$ and $(15)$, we have
$${{1}\over {c^2-\|y\|^2}}\left ( 2\alpha \langle x,y\rangle-|\langle x,y\rangle|^2-{{\alpha^2}\over {c^2}}\|y\|^2\right )<\|x\|^2\ .\eqno{(16)}$$
At this point, taking into account that $c^2\langle \tilde x,f(\tilde x)\rangle\neq \alpha\|f(\tilde x)\|^2$ (by $(b)$), we draw $(12)$ from $(16)$ since
$\tilde x\in I(f(\tilde x),\langle\tilde x,f(\tilde x)\rangle)$.
Summarizing: taking $I={\bf R}$ and $I_n=[-n,n]$ ($n\in {\bf N}$), the continuous function $\varphi$ satisfies $(a_1)$ and $(b_1)$ of Theorem 1,
but, in view of $(7)-(11)$, it does not satisfy the conclusion of that theorem. As a consequence, 
there exists $\tilde\lambda\in {\bf R}$ such that the set of all global minima of
$\varphi(\cdot,\tilde\lambda)$ is disconnected. But such a set agrees with the set of all solutions of the
equation $x=\tilde\lambda f(x)$, and the proof is complete.\hfill $\bigtriangleup$\par
\medskip
REMARK 2. - We do not know whether Theorem 3 is still true when $0\in K$ and $(b)$ is (necessarily) changed in
$$\{x\in X : f(x)\neq 0,\hskip 3pt  c^2\langle x,f(x)\rangle=\alpha\|f(x)\|^2\}\subset \{x\in X : \|f(x)\|\geq c\}\ .$$
However, the proof of Theorem 3 shows that the following is true:\par
\medskip
THEOREM 4. - {\it Let $(X,\langle\cdot,\cdot\rangle)$ be a finite-dimensional real Hilbert space and
let $f:X\to X$ be a continuous function with bounded range.
Assume that there are two numbers $\alpha, c$, with 
$$\inf_{x\in X}\|f(x)\|<c<\|f(0)\|\ ,$$
such that:\par
\noindent
$(a')$\hskip 5pt $\{x\in X : \langle x,f(x)\rangle=\alpha\}\subset \{x\in X : \|f(x)\|<c\}$\ ;\par
\noindent
$(b')$\hskip 5pt $\{x\in X : f(x)\neq 0,\hskip 3pt c^2\langle x,f(x)\rangle=\alpha\|f(x)\|^2\}\subset \{x\in X : \|f(x)\|\geq c\}\ .$\par
Then, there exists $\tilde\lambda\in {\bf R}$ such that the set
$$\{x\in X : x=\tilde\lambda f(x)\}$$
is disconnected.}\par
\medskip
Finally, we present two applications of Theorem 2.\par
\medskip
THEOREM 5. - {\it Let $X$ be a Banach space, let $\varphi\in X^*\setminus \{0\}$ and let $\psi:X\to {\bf R}$ be a Lipschitzian functional
whose Lipschitz constant is equal to $\|\varphi\|_{X^*}$. Moreover, let $[a,b]$ be a compact real interval, $\gamma:[a,b]\to [-1,1]$ a
convex (resp. concave) and continuous function, with $\hbox {\rm int}(\gamma^{-1}(\{-1,1\}))=\emptyset$, and $c\in {\bf R}$. Assume that 
$$\gamma(a)\psi(x)+ca\neq \gamma(b)\psi(x)+cb$$
for all $x\in X$ such that $\psi(x)>0$ (resp. $\psi(x)<0$).\par
Then (with the convention $\sup\emptyset=-\infty$), one has
$$\sup_{\lambda\in \gamma^{-1}(\{-1,1\})}\inf_{x\in X}(\varphi(x)+\gamma(\lambda)\psi(x)+c\lambda)=
\inf_{x\in X}\sup_{\lambda\in [a,b]}(\varphi(x)+\gamma(\lambda)\psi(x)+c\lambda)\ .$$}\par
\smallskip
PROOF. Consider the continuous function $f:X\times [a,b]\to {\bf R}$ defined by
$$f(x,\lambda)=\varphi(x)+\gamma(\lambda)\psi(x)+c\lambda$$
for all $(x,\lambda)\in X\times [a,b]$. By Theorem 2 of [4], for each $\lambda\in \gamma^{-1}(]-1,1[)$, the function
$f(\cdot,\lambda)$ is inf-connected and unbounded below. Also, notice that $\gamma^{-1}(]-1,1[)$, by assumption, is dense in $[a,b]$. Now
fix $x\in X$. If $\psi(x)>0$ (resp. $\psi(x)<0$) the function $f(x,\cdot)$ is convex and, by assumption, $f(x,a)\neq f(x,b)$.
As a consequence, the unique global maximum of this function is either $a$ or $b$. If $\psi(x)\leq 0$, the function is concave
and so, obviously, the set of all its global maxima is connected. Now, the conclusion follows directly from Theorem 2.\hfill
$\bigtriangleup$\par
\medskip
Let $(T,{\cal F},\mu)$ be a $\sigma$-finite measure space, 
 $E$  a real Banach space and $p\geq 1$.
\par
\smallskip
As usual, $L^{p}(T,E)$  denotes the space of all (equivalence
classes of) strongly $\mu$-measurable functions $u : T\rightarrow E$ 
such that
$\int_{T}\parallel u(t)\parallel^{p} d\mu<+\infty$, equipped with
the norm $$\parallel u\parallel_{L^{p}(T,E)}=
\left ( \int_{T}\parallel u(t)\parallel^{p}d\mu\right ) ^{1\over p}\ .$$
\smallskip
A set $D\subseteq L^{p}(T,E)$ is said to be decomposable if, for
every $u,v\in D$ and every $A\in {\cal F}$, the function
 $$t\to
\chi_{A}(t)u(t)+(1-\chi_{A}(t))v(t)$$ belongs to $D$, where $\chi_{A}$
denotes the characteristic function of $A$.\par
\smallskip
A real-valued function on $T\times E$ is said to be a Carat\'eodory function if
it is measurable in $T$ and continuous in $E$.\par
\medskip
THEOREM 6. - {\it Let $(T,{\cal F},\mu)$ be a $\sigma$-finite non-atomic measure space, 
$E$  a real Banach space, $p\in [1,+\infty[$, $X\subseteq L^p(T,E)$ a decomposable set,
$[a,b]$ a compact real interval, $\gamma:[a,b]\to {\bf R}$
a convex (resp. concave) and continuous function.
Moreover, let $\varphi, \psi, \omega:T\times E\to {\bf R}$ be three Carath\'eodory functions such that, for some
$M\in L^1(T)$, $k\in {\bf R}$, one has
$$\max\{|\varphi(t,x)|, |\psi(t,x)|, |\omega(t,x)|\}\leq M(t)+k\|x\|^p$$
for all $(t,x)\in T\times E$ and
$$\gamma(a)\int_T\psi(t,u(t))d\mu+a\int_T\omega(t,u(t))d\mu\neq \gamma(b)\int_T\psi(t,u(t))d\mu+b\int_T\omega(t,u(t))d\mu$$
for all $u\in X$ such that $\int_T\psi(t,u(t))d\mu>0$ (resp. $\int_T\psi(t,u(t))d\mu<0$).\par
Then, one has
$$\sup_{\lambda\in [a,b]}\inf_{u\in X}\left ( \int_T(\varphi(t,u(t))+\gamma(\lambda)\psi(t,u(t)))+\lambda\omega(t,u(t)))d\mu
\right )=$$
$$\inf_{u\in X}\sup_{\lambda\in [a,b]}\left ( \int_T(\varphi(t,u(t))+\gamma(\lambda)\psi(t,u(t)))+\lambda\omega(t,u(t))d\mu
\right )\ .$$}\par
\smallskip
PROOF. The proof goes on exactly as that of Theorem 5. So, one considers the function $f:X\times [a,b]\to {\bf R}$
defined by
$$f(u,\lambda)=\int_T(\varphi(t,u(t))+\gamma(\lambda)\psi(t,u(t)))+\lambda\omega(t,u(t)))d\mu$$
for all $(u,\lambda)\in X\times [a,b]$, and realizes that it satisfies the hypotheses of Theorem 2. In particular, for each $\lambda\in [a,b]$,
the inf-connectedness of the function $f(\cdot,\lambda)$ is due to [6], Th\'eor\`eme 7.\hfill $\bigtriangleup$
\bigskip
\bigskip
{\bf Acknowledgement.} The author has been supported by the Gruppo Nazionale per l'Analisi Matematica, la Probabilit\`a e 
le loro Applicazioni (GNAMPA) of the Istituto Nazionale di Alta Matematica (INdAM) and by the Universit\`a degli Studi di Catania, ``Piano della Ricerca 2016/2018 Linea di intervento 2".\par
\bigskip
\bigskip
\bigskip
\bigskip
\centerline {\bf References}\par
\bigskip
\bigskip
\noindent
[1]\hskip 5pt E. KLEIN and A. C. THOMPSON, {\it Theory of correspondences}, John Wiley $\&$ Sons, 1984.
\smallskip
\noindent
[2]\hskip 5pt B. RICCERI, {\it Some topological mini-max theorems via
an alternative principle for multifunctions}, Arch. Math. (Basel),
{\bf 60} (1993), 367-377.\par
\smallskip
\noindent
[3]\hskip 5pt B. RICCERI, 
{\it Nonlinear eigenvalue problems},  
in ``Handbook of Nonconvex Analysis and Applications'' 
D. Y. Gao and D. Motreanu eds., 543-595, International Press, 2010.\par
\smallskip
\noindent
[4]\hskip 5pt B. RICCERI, {\it On the infimum of certain functionals}, in ``Essays in Mathematics and its Applications -
In Honor of Vladimir Arnold", Th. M. Rassias and P. M. Pardalos eds., 361-367, Springer, 2016.\par
\smallskip
\noindent
[5]\hskip 5pt B. RICCERI, {\it On a minimax theorem: an improvement, a new proof and an overview of its applications},
Minimax Theory Appl., {\bf 2} (2017), 99-152.\par
\smallskip
\noindent
[6]\hskip 5pt J. SAINT RAYMOND, {\it Connexit\'e des sous-niveaux des fonctionnelles int\'egrales}, Rend. Circ. Mat. Palermo,
{\bf 44} (1995), 162-168.\par
\smallskip
\noindent
[7]\hskip 5pt M. SION,  {\it On general minimax theorems}, Pacific J. Math., {\bf 8} (1958), 171-176.\par
\bigskip
\bigskip
\bigskip
\bigskip
Department of Mathematics and Informatics\par
University of Catania\par
Viale A. Doria 6\par
95125 Catania, Italy\par
{\it e-mail address}: ricceri@dmi.unict.it

\bye